\documentclass[10pt]{article}
\usepackage{fullpage}
\usepackage{amsfonts}

\renewcommand{\vec}[1]{{\mathbf #1}}
\newcommand{\vecx}{{\vec x}}

\begin{document}

\hbadness=10000
\vbadness=10000

\newtheorem{theorem}{Theorem}
\newtheorem{corollary}[theorem]{Corollary}
\newtheorem{lemma}[theorem]{Lemma}
\newtheorem{observation}[theorem]{Observation}
\newtheorem{proposition}[theorem]{Proposition}
\newtheorem{definition}[theorem]{Definition}
\newtheorem{claim}[theorem]{Claim}
\newtheorem{fact}[theorem]{Fact}
\newtheorem{assumption}[theorem]{Assumption}

\newcommand{\qed}{\rule{7pt}{7pt}}
\newcommand{\dis}{\mathop{\mbox{\rm d}}\nolimits}
\newcommand{\per}{\mathop{\mbox{\rm per}}\nolimits}
\newcommand{\area}{\mathop{\mbox{\rm area}}\nolimits}
\newcommand{\cw}{\mathop{\rm cw}\nolimits}
\newcommand{\ccw}{\mathop{\rm ccw}\nolimits}
\newcommand{\DIST}{\mathop{\mbox{\rm DIST}}\nolimits}
\newcommand{\OP}{\mathop{\mbox{\it OP}}\nolimits}
\newcommand{\OPprime}{\mathop{\mbox{\it OP}^{\,\prime}}\nolimits}
\newcommand{\ihat}{\hat{\imath}}
\newcommand{\jhat}{\hat{\jmath}}
\newcommand{\abs}[1]{\mathify{\left| #1 \right|}}

\newenvironment{proof}{\noindent{\bf Proof}\hspace*{1em}}{\qed\bigskip}
\newenvironment{proof-sketch}{\noindent{\bf Sketch of Proof}\hspace*{1em}}{\qed\bigskip}
\newenvironment{proof-idea}{\noindent{\bf Proof Idea}\hspace*{1em}}{\qed\bigskip}
\newenvironment{proof-of-lemma}[1]{\noindent{\bf Proof of Lemma #1}\hspace*{1em}}{\qed\bigskip}
\newenvironment{proof-attempt}{\noindent{\bf Proof Attempt}\hspace*{1em}}{\qed\bigskip}
\newenvironment{proofof}[1]{\noindent{\bf Proof
of #1:}}{\qed\bigskip}
\newenvironment{remark}{\noindent{\bf Remark}\hspace*{1em}}{\bigskip}


\newcommand{\FOR}{{\bf for}}
\newcommand{\TO}{{\bf to}}
\newcommand{\DO}{{\bf do}}
\newcommand{\WHILE}{{\bf while}}
\newcommand{\AND}{{\bf and}}
\newcommand{\IF}{{\bf if}}
\newcommand{\THEN}{{\bf then}}
\newcommand{\ELSE}{{\bf else}}

\makeatletter
\def\fnum@figure{{\bf Figure \thefigure}}
\def\fnum@table{{\bf Table \thetable}}
\long\def\@mycaption#1[#2]#3{\addcontentsline{\csname
  ext@#1\endcsname}{#1}{\protect\numberline{\csname
  the#1\endcsname}{\ignorespaces #2}}\par
  \begingroup
    \@parboxrestore
    \small
    \@makecaption{\csname fnum@#1\endcsname}{\ignorespaces #3}\par
  \endgroup}
\def\mycaption{\refstepcounter\@captype \@dblarg{\@mycaption\@captype}}
\makeatother

\newcommand{\figcaption}[1]{\mycaption[]{#1}}
\newcommand{\tabcaption}[1]{\mycaption[]{#1}}
\newcommand{\head}[1]{\chapter[Lecture \##1]{}}
\newcommand{\mathify}[1]{\ifmmode{#1}\else\mbox{$#1$}\fi}
\newcommand{\bigO}O
\newcommand{\set}[1]{\mathify{\left\{ #1 \right\}}}
\def\half{\frac{1}{2}}


\newcommand{\enc}{{\sf Enc}}
\newcommand{\dec}{{\sf Dec}}
\newcommand{\E}{{\rm Exp}}
\newcommand{\Var}{{\rm Var}}
\newcommand{\Z}{{\mathbb Z}}
\newcommand{\F}{{\mathbb F}}
\newcommand{\integers}{{\mathbb Z}^{\geq 0}}
\newcommand{\R}{{\mathbb R}}
\newcommand{\Q}{{\cal Q}}
\newcommand{\eqdef}{{\stackrel{\rm def}{=}}}
\newcommand{\from}{{\leftarrow}}
\newcommand{\vol}{{\rm Vol}}
\newcommand{\poly}{{\rm poly}}
\newcommand{\ip}[1]{{\langle #1 \rangle}}
\newcommand{\wt}{{\rm wt}}
\renewcommand{\vec}[1]{{\mathbf #1}}
\newcommand{\mspan}{{\rm span}}
\newcommand{\rs}{{\rm RS}}
\newcommand{\RM}{{\rm RM}}
\newcommand{\Had}{{\rm Had}}
\newcommand{\calc}{{\cal C}}
\newcommand{\binom}[2]{{#1 \choose #2}}

\newcommand{\fig}[4]{
        \begin{figure}
        \setlength{\epsfysize}{#2}
        \vspace{3mm}
        \centerline{\epsfbox{#4}}
        \caption{#3} \label{#1}
        \end{figure}
        }

\newcommand{\ord}{{\rm ord}}

\providecommand{\norm}[1]{\lVert #1 \rVert}
\newcommand{\embed}{{\rm Embed}}
\newcommand{\qembed}{\mbox{$q$-Embed}}
\newcommand{\calh}{{\cal H}}
\newcommand{\lp}{{\rm LP}}
\title{Improved lower bound on the size of Kakeya sets over finite fields}

\author{
Shubhangi Saraf\thanks{MIT CSAIL. {\tt shibs@mit.edu}.}
\and
Madhu Sudan\thanks{MIT CSAIL. {\tt madhu@mit.edu}.}
}

\maketitle


\begin{abstract}

In a recent breakthrough, Dvir showed that every Kakeya set in
$\F^n$ must be of cardinality at least $c_n |\F|^n$ where $c_n \approx
1/n!$. We improve this lower bound to
$\beta^n |\F|^n$ for a constant $\beta > 0$.
This pins down the growth of the leading constant to the
right form as a function of $n$.
\end{abstract}

Let $\mathbb F$ be a finite field of $q$ elements.

\begin{definition}[Kakeya Set]
A set $K \subseteq \F^n$ is said to be a {\em Kakeya set} in $\F^n$,
if for
every $\vec b \in \F^n$, there exists a point $\vec a \in \F^n$
such that for every $t \in \F$, the point $\vec a + t \cdot \vec b \in K$.
\end{definition}

We show:
\begin{theorem}
\label{thm}
There exist constants $c_0,c_1 > 0$ such that for all
$n$, if $K$ is a Kakeya set in $\F^n$ then
$|K| \geq c_0 \cdot (c_1 \cdot q)^n$.
\end{theorem}

\begin{remark}
Our proofs give some tradeoffs on the constants $c_0, c_1$
that are achievable. We comment on the constants at the end of
the paper.
\end{remark}

The question of establishing lower bounds on the size of Kakeya sets
was posed in Wolff~\cite{Wolff}. Till recently, the best known lower
bound on the size of Kakeya sets was of the form $q^{\alpha n}$ for
some $\alpha < 1$. In a recent breakthrough Dvir~\cite{Dvir} showed
that every Kakeya set must have cardinality at least $c_n q^n$ for
$c_n = (n!)^{-1}$. (Dvir's original bound achieved a weaker lower
bound of $c_n \cdot q^{n-1}$, but \cite{Dvir} includes the stronger
bound of $c_n \cdot q^n$, with the improvements being attributed to
Alon and Tao.)
Our improvement now shows that $c_n$ remains bounded from below by
$\beta^n$ for some fixed $\beta > 0$.
While this improvement in the lower bound on size of Kakeya sets
is quantitatively small (say, compared to
the improvement of Alon and Tao over Dvir's original bound) it
is qualitatively significant in that it
does determine the growth of the leading constant $c_n$, up to the
determination of the right constant $\beta$.
In particular, it compares well with known upper bounds.
Previously, it was known there exists a constant $\beta
< 1$ such that there are Kakeya sets of cardinality at most $\beta^n
q^n$, for every odd $q$.
A bound of
$\beta \leq \frac{1}{\sqrt{2}}$ follows from
Mockenhaupt and Tao~\cite{MoTa}
and the fact that products of Kakeya sets are Kakeya sets
(in higher dimension).
The best known constant has $\beta \to \frac12$
due to Dvir~\cite{Dvir:personal-comm}.
We include his proof in the appendix (see Appendix~\ref{app}),
complementing it with a similar construction and bound for
the case of even $q$ as well (so now the upper bounds work for all
large fields).

Our proof follows that of Dvir~\cite{Dvir}.
Given a Kakeya set $K$ in $\F^n$, we
show that there exists an $n$-variate polynomial, whose degree is bounded
from above by some function of $|K|$, that
vanishes at all of $K$. Looking at restrictions of this polynomial
on lines yields that this polynomial has too many zeroes, which in turn
yields a lower bound on the size of $K$.
Our main difference is that we look for polynomials that vanish with
``high multiplicity'' at each point in $K$. The requirement of high
multiplicity forces the degree of the $n$-variate polynomial to go up
slightly, but yields more zeroes when this polynomial is restricted to
lines. The resulting tradeoff turns out to yield an improved bound.
(We note that this is
similar to
the techniques used for ``improved list-decoding of Reed-Solomon codes''
by Guruswami and Sudan~\cite{GuSu}.)
We now give a formal proof of Theorem~\ref{thm}.

\section{Preliminaries}

For $\vecx = \langle x_1,\ldots,x_n \rangle$, let $\F[\vecx]$
denote the ring of polynomials in $x_1,\ldots,x_n$ with coefficients
from $\F$. We recall the following basic fact on polynomials.

\begin{fact}
\label{fact}
Let $P \in \F[\vecx]$ be a polynomial of degree at most $q-1$
in each variable. If $P(\vec a) = 0$ for all $\vec a \in \F^n$,
then $P \equiv 0$.
\end{fact}

For real $\alpha \geq 0$, let $N_q(n,m)$ denote the number of
monomials in $n$ variables of total degree less than
$m \cdot q$ and of individual degree at most $q-1$ in
each variable.

We say that a polynomial $g \in \F[\vecx]$ has a zero of
{\em multiplicity} $m$ at a
point $\vec a \in \F^n$ if the polynomial $g_{\vec a}(\vec x) =
g(\vecx + \vec a)$ has no support on
monomials of
degree strictly less than $m$.
Note that the coefficients of $g_{\vec a}$ are (homogenous) linear
forms in
the coefficients of $g$ and thus the constraint $g$ has a zero of
multiplicity $m$ at $\vec a$ yields ${m + n -1 \choose n}$ homogenous
linear
constraints on the coefficients of $g$. As a result we conclude:

\begin{proposition}
\label{propexists}
Given a set $S \subseteq \F^n$ satisfying ${m + n-1 \choose n} \cdot
|S| < N_q(n,m)$, there exists a non-zero polynomial
$g \in \F[\vecx]$ of total degree less than $mq$ and degree
at most $q-1$ in each variable such that $g$ has a zero of
multiplicity $m$ at every point $\vec a \in S$.
\end{proposition}

\begin{proof}
The number of possible coefficients for $g$ is $N_q(n,m)$ and
the number of (homogenous linear) constraints is
${m + n-1 \choose n} \cdot |S| < N_q(n,m)$.
Since the number of constraints is strictly smaller than the
number of unknowns, there is a non-trivial solution.
\end{proof}

For $g \in \F[\vecx]$ we let restriction $g_{\vec a, \vec b}(t)
= g(\vec a + t\cdot \vec b)$
denote its restriction to the ``line'' $\{\vec a + t\cdot \vec b | t\in
\F\}$. We note the following facts on the restrictions of polynomials
to lines.

\begin{proposition}
\label{propmult}
If $g \in \F[\vec x]$ has a root of multiplicity $m$ at
some point $\vec a + t_0 \vec b$ then $g_{\vec a, \vec b}$ has
a root of multiplicity $m$ at $t_0$.
\end{proposition}

\begin{proof}
By definition, the fact that $g$ has a zero of multiplicity $m$ at
$\vec a + t_0 \vec b$ implies that
the polynomial $g(\vecx + (\vec a + t_0\vec b))$
has no support on monomials of degree less than $m$.
Thus, under the homogenous substitution $\vec x \leftarrow t\cdot \vec b$,
we get no monomials of degree less than $m$ either, and thus we have
$t^m$ divides $g(t\vec b + (\vec a + t_0 \vec b)) =
g(\vec a + (t+t_0) \vec b) = g_{\vec a, \vec b}(t + t_0)$.
The final form implies that that $g_{\vec a, \vec b}$ has a zero
of multiplicity $m$ at $t_0$.
\end{proof}

\begin{proposition}[\cite{Dvir}]
\label{proplead}
Let $g \in \F[\vec x]$ be a non-zero polynomial of total
degree $d$ and let $g_0$ be the (unique, non-zero) homogenous
polynomial
of degree $d$ such that $g = g_0 + g_1$ for some polynomial
$g_1$ of degree strictly less than $d$.
Then $g_{\vec a, \vec b}(t) = g_0(\vec b) t^d + h(t)$ where
$h$ is a polynomial of degree strictly less than $d$.
\end{proposition}

\section{Proof of Theorem~\ref{thm}}

\begin{lemma}
\label{lemmain}
If $K$ is a Kakeya set in $\F^n$, then
for every integer $m \geq 0$,
$|K| \geq \frac{1}{\binom{m+n-1}{n}} \cdot N_q(n,m)$.
\end{lemma}
\begin{proof}
Assume for contradiction that
$|K| < \frac{1}{\binom{m+n-1}{n}} \cdot N_q(n,m)$.
Let $g \in \F[\vecx]$ be a non-zero polynomial of total degree
less than $mq$ and degree at
most $q-1$ in each variable
that has a zero of multiplicity $m$ for each $x \in K$. (Such
a polynomial exists by Proposition~\ref{propexists}.)
Let $d < mq$ denote the total degree of $g$ and let $g =
g_0 + g_1$ where $g_0$ is homogenous of degree
$d$ and $g_1$ has degree less than $d$. Note that $g_0$ is
also non-zero and has degree at most $q-1$ in every variable.

Now fix a ``direction'' $\vec b \in \F^n$. Since $K$ is a Kakeya
set, there exists $\vec a \in \F^n$ such that $\vec a + t \vec b \in
K$ for every $t \in \F$. Now consider the restriction $g_{\vec
a,\vec b}$ of $g$ to the line through $\vec a$ in direction $\vec
b$. $g_{\vec a,\vec b}$ is a univariate polynomial of degree at most
$d < mq$. At every point $t_0 \in \F$ we have that $g_{\vec a,\vec b}$
has a zero of multiplicity $m$ (Proposition~\ref{propmult}). Thus
counting up the zeroes of $g_{\vec a, \vec b}$ we find it has $mq$
zeroes ($m$ at every $t_0 \in \F$) which is more than its degree.
Thus $g_{\vec a,\vec b}$ must be identically zero. In particular
its leading coefficient must be zero. By
Proposition~\ref{proplead} this leading coefficient is $g_0(\vec b)$
and so we conclude $g_0(\vec b) = 0$.

We conclude that $g_0$ is zero on all of $\F^n$ which contradicts
the fact (Fact~\ref{fact}) that it is a non-zero polynomial of
degree at most $q-1$ in each of its variables.
\end{proof}

\begin{proofof}{Theorem~\ref{thm}}
Theorem~\ref{thm} now follows by choosing $m$ appropriately.
Using for instance $m = n$, we obtain that
$|K| \geq \frac{1}{\binom{2n-1}n} \cdot q^n \geq (q/4)^n$,
establishing the theorem for $c_0 = 1$ and $c_1 = 1/4$.

A better choice is with $m = \lceil n/2 \rceil \leq (n+1)/2$.
In this case $N_q(n,m) \geq
\frac12 q^n$ (since at least half the monomials of individual degree at
most $q-1$ have degree at most $nq/2$). This leads to a bound
of $|K| \geq \frac{1}{2\binom{(3/2)n}n} q^n \geq \frac12 (q/2.6)^n$,
yielding the theorem for $c_0 = 1/2$ and $c_1 = 1/2.6$.
\end{proofof}

To improve the constant $c_1$ further, one could
study the asymptotics of $N_q(n,m)$ closer. Let
$\tau_{\alpha}$ denote the quantity
$\liminf_{n \to \infty} \{\liminf_{q \to \infty} \frac1q\cdot
N_q(n,\alpha\cdot n)^{1/n}\}$.
I.e., for sufficiently large $n$ and sufficiently larger $q$,
$N_q(n,\alpha n) \to \tau_{\alpha}^n \cdot q^n$.
Lemma~\ref{lemmain} can be reinterpreted in these terms as saying
that for every $\alpha \in [0,1]$, every Kakeya set has size
at least $c_0 (c_\alpha \cdot q)^n - o(q^n)$
for some $c_0 > 0$,
where $c_\alpha \to \tau_{\alpha}/2^{(1+\alpha)H(1/(1+\alpha)}$
(where $H(x) = - x \log_2 x - (1-x)\log_2 (1-x)$ is the binary
entropy function).
The best estimate on $\tau_{\alpha}$ we were able to obtain
does not have a simple closed form expression. As $q \to \infty$,
$\tau_{\alpha}^n$ equals the
volume of the following region in $\R^n$: $\{(x_1,x_2,\ldots x_n) \in [0,1]^n| \sum_{i=1}^n x_i \leq \alpha \cdot n\}$. This volume can be expressed in terms of Eulerian numbers (See~\cite{JeMo}, $\S 4.3$).
\cite[$\S 6$]{GiKe} gives some asymptotics for Eulerian numbers
and using their estimates $\alpha = 0.398$, it seems one can
reduce $c_{\alpha}$ to something like $\frac{1}{2.46}$. This still
remains bounded away from the best known upper bound which has
$c_1 \to 1/2$.

\begin{remark}
While the main theorem only gives the limiting behavior of
Kakeya sets for large $n,q$, Lemma~\ref{lemmain} can still be
applied to specific choices and get improvements over~\cite{Dvir}.
For example, for $n=3$, using
$m=2$ we get a lower bound of $\frac{5}{24} q^3$ as opposed to
the bound of $\frac16 q^3$ obtainable from \cite{Dvir}.
\end{remark}

\section*{Acknowledgments}

Thanks to Zeev Dvir for explaining the Kakeya problem and his solution
to us, for detailed answers to many queries, and for his permission
to include his upper bound on the size of Kakeya sets here (see Appendix~\ref{app}).
Thanks also to Swastik Kopparty for helping us extend Dvir's proof to even
characteristic.
Thanks to Chris Umans and Terry Tao for valuable discussions.

\bibliographystyle{plain}

\begin{thebibliography}{1}

\bibitem{Dvir}
Zeev Dvir.
\newblock On the size of {K}akeya sets in finite fields.
\newblock {\em Journal of the American Mathematical Society}, (to appear),
  2008.
\newblock Article electronically published on June 23, 2008.

\bibitem{Dvir:personal-comm}
Zeev Dvir.
\newblock Personal communication, August 2008.

\bibitem{GiKe}
Eldar Giladi and Joseph~B. Keller.
\newblock Eulerian number asymptotics.
\newblock {\em Proceedings of the Royal Society of London, Series A},
  445(1924):291--303, 1994.

\bibitem{GuSu}
Venkatesan Guruswami and Madhu Sudan.
\newblock Improved decoding of {R}eed-{S}olomon and algebraic-geometric codes.
\newblock {\em IEEE Transactions on Information Theory}, 45:1757--1767, 1999.

\bibitem{JeMo}
Jean-Luc Marichal and Michael~J. Mossinghoff.
\newblock Slices, slabs, and sections of the unit hypercube.
\newblock {\em Online Journal of Analytic Combinatorics}, 3, 2008.
\newblock 11 pages. Earlier version appears as eprint math/0607715 (2006).

\bibitem{MoTa}
Gerd Mockenhaupt and Terence Tao.
\newblock Restriction and {K}akeya phenomena for finite fields.
\newblock {\em Duke Mathematics Journal}, 121(1):35--74, 2004.

\bibitem{Wolff}
T.~Wolff.
\newblock Recent work connected with the {K}akeya problem.
\newblock In {\em Prospects in Mathematics}, pages 129--162. Princeton, NJ,
  1999.

\end{thebibliography}

\appendix

\section{An upper bound on Kakeya sets}
\label{app}

We include here Dvir's proof~\cite{Dvir:personal-comm} giving
a non-trivial upper bound on the size of Kakeya sets in fields
of odd characteristic. The proof is based on the construction
of Mockenhaupt and Tao~\cite{MoTa}.
For the case of even characteristic we complement their results
by using a variation (obtained with Swastik Kopparty) of their
construction.

\begin{theorem}[\cite{Dvir:personal-comm}]
For every $n\geq 2$, and field $\F$, there
exists a Kakeya set in $\F^n$ of cardinality at most
$2^{-(n-1)}\cdot q^n + O(q^{n-1})$.
\end{theorem}

\begin{proof}
We consider two cases depending on whether $\F$ is of odd or
even characteristic.

\noindent\textbf{Odd characteristic:}
Let $D_n = \{\langle \alpha_1,\ldots,\alpha_{n-1},\beta\rangle
| \alpha_i,\beta \in \F, \alpha_i + \beta^2$ is a square$\}$.
Now let $K_n = D_n \cup (\F^{n-1} \times \{0\})$ where
$\F^{n-1} \times \{0\}$ denotes the set $\{\langle \vec a,0\rangle
| \vec a \in \F^{n-1}\}$.
We claim that $K_n$ is a Kakeya set of the appropriate size.

Consider a direction $\vec b = \langle b_1,\ldots,b_n\rangle$.
If $b_n = 0$, for $\vec a = \langle 0,\ldots,0\rangle$ we
have that $\vec a + t \vec b \in \F^{n-1} \times \{0\}
\subseteq K_n$.
The more interesting case is when $b_n \ne 0$. In this case
let $\vec a = \langle (b_1/(2b_n))^2,\ldots,(b_{n-1}/(2b_n))^2,0
\rangle$. The point $\vec a + t \vec b$ has coordinates
$\langle \alpha_1,\ldots,\alpha_{n-1},\beta\rangle$
where $\alpha_i = (b_i/(2b_n))^2 + tb_i$ and
$\beta = tb_n$. We have $\alpha_i + \beta^2 = (b_i/(2b_n) + tb_n)^2$
which is a square for every $i$ and so
$\vec a + t \vec b \in D_n \subseteq K_n$.
This proves that $K_n$ is indeed a Kakeya set.

Finally we verify that the size of $K_n$ is as claimed.
First note that the size of $D_n$ is exactly
$q \cdot ((q+1)/2)^{n-1} = 2^{-(n-1)} q^n + O(q^{n-1})$
($q$ choices for $\beta$ and $(q+1)/2$ choices for each
$\alpha_i + \beta^2$).
The size of $K_n$ is at most $|D_n| + q^{n-1}
= 2^{-(n-1)} q^n + O(q^{n-1})$ as claimed.

\noindent\textbf{Even characteristic:}
This case is handled similarly with minor variations in the
definition of $K_n$.
Specifically, we let
$K_n = E_n = \{\langle \alpha_1,\ldots,\alpha_{n-1},\beta\rangle
| \alpha_i,\beta \in \F, \exists \gamma_i \in \F$ such that $\alpha_i = \gamma_i^2 + \gamma_i\beta\}$.
(As we see below $E_n$ contains $\F^{n-1} \times \{0\}$ and so
there is no need to set $K_n = E_n \cup \F^{n-1} \times \{0\}$.)

Now consider direction $\vec b = \langle b_1,\ldots,b_n\rangle$.
If $b_n = 0$, then let $\vec a = 0$. We note that
$\vec a + t\vec b = \langle t b_1, \ldots, t b_{n-1}, 0\rangle
=
\langle \gamma_1^2 + \beta \gamma_1,\ldots
\gamma_{n-1}^2 + \beta \gamma_{n-1},\beta\rangle$
for $\beta = 0$ and $\gamma_i = \sqrt{tb_i} = (tb_i)^{q/2}$.
We conclude that
$\vec a + t\vec b \in E_n$ for every $t \in \F$ in this case.
Now consider the case where $b_n \neq 0$.
Let $\vec a = \langle (b_1/b_n)^2,\ldots,(b_{n-1}/b_n)^2,0 \rangle$.
The point $\vec a + t \vec b$ has coordinates
$\langle \alpha_1,\ldots,\alpha_{n-1},\beta\rangle$
where $\alpha_i = (b_i/b_n)^2 + tb_i$ and
$\beta = tb_n$.  For $\gamma_i = (b_i/b_n)$,
$\gamma_i^2 + \gamma_i\beta = (b_1/b_n)^2 + tb_i = \alpha_i$.
Hence $\vec a + t \vec b \in E_n = K_n$.

It remains to compute the size of $E_n$.
The number of points of the form
$\langle \alpha_1,\ldots,\alpha_{n-1},0\} \in E_n$
is exactly $q^{n-1}$.
We now determine the size of
$\langle \alpha_1,\ldots,\alpha_{n-1},\beta\} \in E_n$ for
fixed $\beta \ne 0$.
We first claim that the set $\{\gamma^2 + \beta \gamma | \gamma \in \F\}$
has size exactly $q/2$. This is so since for every
$\gamma \in \F$, we have
$\gamma^2 + \beta \gamma = \tau^2 + \beta \tau$ for
$\tau = \gamma + \beta \ne \gamma$, and so the map
$\gamma \mapsto \gamma^2 + \beta \gamma$ is a 2-to-1 map on
its image.
Thus, for $\beta \ne 0$,
the number of points of the form
$\langle \alpha_1,\ldots,\alpha_{n-1},\beta\}$ in $E_n$ is
exactly $(q/2)^{n-1}$.
We conclude that $E_n$ has cardinality $(q-1) \cdot (q/2)^{n-1} +
q^{n-1} = 2^{-(n-1)} q^n + O(q^{n-1})$.
\end{proof}

We remark that for the case of odd characteristic, one can also
use a recursive construction, replacing the set
$\F^{n-1}\times \{0\}$ by $K_{n-1} \times \{0\}$. This would reduce
the constant in the $O(q^{n-1})$ term, but not alter the leading
term. Also we note that the construction used in the even case
essentially also works in the odd characteristic case.
Specifically the set $E_n \cup \F^{n-1} \times \{0\}$ is a
Kakeya set also for odd characteristic. Its size can also be
argued to be $2^{-(n-1)}\cdot q^n + O(q^{n-1})$.

\end{document}